\documentclass[10pt]{article}
\usepackage[T1]{fontenc}
\usepackage[utf8]{inputenc}
\usepackage{lmodern}
\usepackage[a4paper,margin=1in]{geometry}
\usepackage{amsmath,amssymb,amsthm,mathtools,mathrsfs}
\usepackage{microtype}
\usepackage[skip=0.45\baselineskip plus 2pt,indent=0pt]{parskip}
\usepackage{enumitem}
\usepackage[colorlinks=true,linkcolor=blue,citecolor=blue,urlcolor=blue]{hyperref}
\usepackage{aliascnt}
\usepackage[nameinlink,noabbrev,capitalise]{cleveref}
\allowdisplaybreaks[2]
\newtheorem{theorem}{Theorem}[section]
\newaliascnt{fact}{theorem}
\newtheorem{fact}[fact]{Fact}
\aliascntresetthe{fact}
\crefname{fact}{Fact}{Facts}
\Crefname{fact}{Fact}{Facts}
\newaliascnt{lemma}{theorem}
\newtheorem{lemma}[lemma]{Lemma}
\aliascntresetthe{lemma}
\crefname{lemma}{Lemma}{Lemmas}
\Crefname{lemma}{Lemma}{Lemmas}
\newaliascnt{claim}{theorem}
\newtheorem{claim}[claim]{Claim}
\aliascntresetthe{claim}
\crefname{claim}{Claim}{Claims}
\Crefname{claim}{Claim}{Claims}
\newaliascnt{definition}{theorem}
\newtheorem{definition}[definition]{Definition}
\aliascntresetthe{definition}
\numberwithin{equation}{section}
\newcommand{\N}{\mathbb N}
\newcommand{\E}{\mathbb E}\newcommand{\Prb}{\mathbb P}\newcommand{\eps}{\varepsilon}
\newcommand{\bx}{\boldsymbol x}\newcommand{\ba}{\boldsymbol a}
\newcommand{\F}{\mathcal F}
\newcommand{\unif}{\mathrm{unif}}

\setlist[enumerate]{label=\textup{(\roman*)},itemsep=2pt,topsep=4pt,leftmargin=2em}
\hypersetup{
 pdftitle={Hypergraphs with zero uniform Turan density at intermediate levels},
 pdfsubject={Zero uniform Turan density at intermediate levels and position-consistent orderings},
 pdfauthor={Heng Li and Xizhi Liu},
 bookmarksnumbered=true}
\title{\Large\bfseries Hypergraphs with zero uniform Tur\'an density at~intermediate~levels}
\author{%
Heng~Li%
\and
Xizhi~Liu%
}

\date{22 October 2026}
\begin{document}
\maketitle

\begin{abstract}
An $r$-uniform hypergraph $H$ is uniformly dense with respect to $j$-graphs if, up to a small error, at least a prescribed proportion of the $r$-cliques of every $j$-uniform hypergraph on the same vertex set are edges of $H$. For a forbidden hypergraph $F$, its uniform Tur\'an density at level $j$ is the supremum of these density bounds over arbitrarily large $F$-free hypergraphs with the error tending to zero. For every $r\ge4$ and $1\le j\le r-3$, this paper proves that this Tur\'an density is zero precisely when $F$ admits a vertex ordering in which every $(j+1)$-set shared by two edges occupies the same positions in both edges. This generalizes the zero-density characterizations of Erd\H{o}s at $j=0$ and of Reiher, R\"odl, and Schacht at $j=r-2$, and establishes the latter authors' proposed extension to the remaining nontrivial levels. Consequently, at level $j$, every positive uniform Tur\'an density of a possibly infinite forbidden family is at least $N^{-N}$, where $N=\binom r{j+1}$. The proof combines hypergraph regularity with adjacent exchanges in spreadable random hypergraphs.
\end{abstract}

\section{Introduction}
For a nonnegative integer $r$, an \emph{$r$-uniform hypergraph}, or \emph{$r$-graph}, $H$ consists of a nonempty finite vertex set $V(H)$ and a collection of $r$-element subsets of $V(H)$, called edges. We identify a hypergraph with its edge set and write $v(H)\coloneqq|V(H)|$ for its number of vertices and $|H|$ for its number of edges. For a positive integer $m$, write $[m]\coloneqq\{1,\ldots,m\}$. For a finite set $V$, let $\binom Vr$ denote the family of its $r$-element subsets. Given integers $r\ge1$ and $0\le j\le r$, and a $j$-graph $G$ on $V$, put
\[
 \mathcal K_r(G)\coloneqq\left\{e\in\tbinom Vr\colon\tbinom ej\subseteq G\right\}.
\]
For $r\ge3$, $0\le j<r$, $d\in[0,1]$, and $\eta>0$, an $r$-graph $H$ on $n$ vertices is \emph{$(d,\eta,j)$-dense} if, for every $j$-graph $G$ on $V(H)$,
\begin{equation}\label{eq:dense}
 |H\cap\mathcal K_r(G)|\ge d|\mathcal K_r(G)|-\eta n^r.
\end{equation}
Intuitively, every auxiliary $j$-graph $G$ selects the $r$-sets whose $j$-subsets all lie in $G$, and at least a proportion $d$ of these selected $r$-sets must be edges of $H$, up to an additive error of $\eta n^r$ in the number of edges.

For a possibly infinite family $\F$ of $r$-graphs, the \emph{uniform Tur\'an density at level $j$} of $\F$ is the supremum of all $d\in[0,1]$ such that, for every $\eta>0$ and $n_0\in\N$, there is an $\F$-free $(d,\eta,j)$-dense $r$-graph on at least $n_0$ vertices; that is,
\[
 \pi_{j,\unif}(\F)\coloneqq\sup\left\{d\in[0,1]\colon
 \begin{aligned}
  &\forall\eta>0\ \forall n_0\in\N\ \exists\,r\text{-graph }H\text{ with }v(H)\ge n_0\\
  &\text{such that }H\text{ is }\F\text{-free and }(d,\eta,j)\text{-dense}
 \end{aligned}
 \right\}.
\]
For a single forbidden hypergraph $F$, write $\pi_{j,\unif}(F)\coloneqq\pi_{j,\unif}(\{F\})$. We assume that every forbidden hypergraph has an edge; the family $\F$ may be empty.

At $j=0$, taking $G=\{\varnothing\}$ reduces \eqref{eq:dense} to $|H|\ge d\binom nr-\eta n^r$, and the empty $0$-graph imposes no further restriction. Thus this level recovers the ordinary Tur\'an problem~\cite{Turan1941} of maximizing the number of edges in an $n$-vertex $r$-graph that contains no copy of a fixed $r$-graph $F$. Writing $\operatorname{ex}(n,F)$ for this maximum, the \emph{ordinary Tur\'an density} of $F$ is $\pi(F)\coloneqq\lim_{n\to\infty}{\operatorname{ex}(n,F)}/{\binom nr}=\pi_{0,\unif}(F)$. For ordinary graphs, the Erd\H{o}s--Stone--Simonovits theorem~\cite{ErdosSimonovits,ErdosStone} determines the Tur\'an density from the chromatic number. In higher uniformities, even the ordinary Tur\'an density of the tetrahedron remains unresolved; see the survey of Keevash~\cite{KeevashSurvey}. For the zero-density question, Erd\H{o}s~\cite{ErdosExtremal} proved that $\pi(F)=0$ if and only if $F$ is $r$-partite. Thus these are precisely the $r$-graphs $F$ such that, for every $d>0$, every sufficiently large $r$-graph $H$ with $|H|\ge d\binom{v(H)}r$ contains a copy of $F$.

At $j=1$, every auxiliary $1$-graph has the form $G=\binom U1$ for some $U\subseteq V(H)$, so $\mathcal K_r(G)=\binom Ur$. The density condition therefore requires $|H\cap\binom Ur|\ge d\binom{|U|}r-\eta n^r$ for every $U\subseteq V(H)$. For every fixed $\alpha\in(0,1]$ and sufficiently large $n$, this gives an induced edge density of at least $d-O_r(\eta/\alpha^r)$ on every set $U$ with $|U|\ge\alpha n$. Erd\H{o}s and S\'os~\cite{ErdosSos} initiated uniform Tur\'an problems in this setting, asking which hypergraphs must occur under a uniform lower bound on the edge densities of large vertex subsets. For $3$-graphs, Glebov, Kr{\'a}{\v{l}}, and Volec~\cite{GKV} proved that the tetrahedron with one edge removed has uniform Tur\'an density $1/4$, and Reiher, R\"odl, and Schacht~\cite{RRSweak} gave another proof. More recently, Kielak et al.~\cite{KKLLSWTetrahedron} proved that the tetrahedron itself has uniform Tur\'an density $1/2$, and Buci\'c~\cite{BucicTetrahedron} gave a Fourier-analytic proof of the same result. For general forbidden $3$-graphs, palette constructions color vertex pairs and select triples according to allowed color patterns. Lamaison~\cite{LamaisonPalettes} proved that $F$-free palette constructions yield densities arbitrarily close to $\pi_{1,\unif}(F)$ from below for every forbidden $3$-graph $F$. Turning to the zero-density question, the class of forbidden $r$-graphs with zero density at level $1$ extends beyond the $r$-partite hypergraphs. Kohayakawa, Nagle, R\"odl, and Schacht~\cite{KNRSlinear} proved that every linear $r$-graph has zero uniform Tur\'an density at level $1$. For $3$-graphs, Ding et al.~\cite{DLLWYCodegree} showed that vanishing codegree Tur\'an density forces vanishing uniform Tur\'an density at level $1$.

For general $j$, Reiher, R\"odl, and Schacht~\cite[Section~1.4]{RRS} introduced the density conditions in \eqref{eq:dense} and the corresponding Tur\'an densities; see also the survey of Reiher~\cite{ReiherSurvey}. For every fixed forbidden $r$-graph $F$, these densities satisfy~\cite[Equation~(1.3)]{RRS}
\[
 \pi_{0,\unif}(F)\ge\pi_{1,\unif}(F)\ge\cdots\ge\pi_{r-2,\unif}(F)\ge\pi_{r-1,\unif}(F)=0.
\]
Increasing $j$ strengthens the density condition on the host, and the last nontrivial level is $j=r-2$. At this level, Reiher, R\"odl, and Schacht~\cite[Theorem~1.8]{RRS} proved that the $r$-graph on $r+1$ vertices with three edges has uniform Tur\'an density $2^{1-r}$. Lin, Sun, Wang, and Zhou~\cite{LSWZ} established a palette characterization for forbidden families at this level and obtained further exact values. Related notions based on vertex-pair incidences, two-paths, and quasirandom links were studied by Reiher, R\"odl, and Schacht~\cite{RRStetrahedra,RRStwopath}, Berger et al.~\cite{BPRRSQuasirandomLinks}, and Lin, Ma, and Zhou~\cite{LMZVanishing}. Returning to the density condition in \eqref{eq:dense}, Reiher, R\"odl, and Schacht~\cite[Theorems~1.2 and~6.1]{RRSzero} gave a complete characterization of zero density at $j=r-2$. An $r$-graph has zero density at this level if and only if its vertices can be ordered so that every $(r-1)$-set occupies the same positions in all edges containing it. Reiher, R\"odl, and Schacht~\cite[Section~6.1]{RRSzero} proposed the analogous characterization for $1\le j\le r-3$, and Lamaison~\cite[Conjecture~20]{LamaisonBeyond} also conjectured the case $j=1$. 

For $2\le q\le r-1$, call a vertex ordering of an $r$-graph \emph{$q$-position-consistent} if every $q$-set contained in two edges occupies the same positions when the vertices of each edge are listed in increasing order. The main result establishes this characterization.

\begin{theorem}\label{thm:zero}
Let $r\ge4$, $1\le j\le r-3$, $q=j+1$, and $N=\binom rq$. An $r$-graph $F$ satisfies $\pi_{j,\unif}(F)=0$ if and only if it has a $q$-position-consistent ordering. Consequently, for every family $\F$ of $r$-graphs, either $\pi_{j,\unif}(\F)=0$ or $\pi_{j,\unif}(\F)\ge N^{-N}$.
\end{theorem}

Together with the classical case $j=0$ and the result at $j=r-2$, this settles the zero-density question at every nontrivial level of the hierarchy above.

At $j=0$, Erd\H{o}s~\cite{ErdosExtremal} established the lower bound $r!/r^r$ for positive densities, and Mubayi~\cite{MubayiTuranExtension} showed that it is attained by a single forbidden $r$-graph. At $j=1$, Lamaison~\cite[Theorem~3]{LamaisonBeyond} constructed a single forbidden $r$-graph whose uniform Tur\'an density at level $1$ is $\binom r2^{-\binom r2}$. At $j=r-2$, the lower bound $r^{-r}$ is attained by a single forbidden $r$-graph, as shown by Garbe, Kr{\'a}{\v{l}}, and Lamaison~\cite{GKLMinimum} for $r=3$ and Lin, Wang, and Zhou~\cite{LWZMinimum} for every $r\ge3$. In a forthcoming survey, we will show that, for all sufficiently large $r$ and every $1\le j\le r-3$, every value in $[N^{-N},1]$ is attained as the uniform Tur\'an density at level $j$ of a possibly infinite forbidden family. In particular, the endpoint $N^{-N}$ is attained, so the lower bound in \cref{thm:zero} is best possible in this range.

\Cref{sec:tools} introduces spreadable random hypergraphs and proves \cref{lem:array-embedding}. The regularity tools and the proof of \cref{thm:zero} are given in \cref{sec:zero}. All hypergraphs are finite and simple unless stated otherwise.

\section{An embedding theorem for spreadable random hypergraphs}\label{sec:tools}

In this section, we consider $q$-shadow-independent spreadable random $r$-graphs, defined below. This random $r$-graph model arises from the regularity reduction in \cref{sec:zero}. The embedding theorem below shows that positive edge probability forces every ordered $q$-position-consistent $r$-graph to occur with positive probability, giving the contradiction needed for the zero-density direction of \cref{thm:zero}.

To state these properties precisely, we represent a random $r$-graph on an ordered vertex set $V$ by its edge indicators $(X_e)_{e\in\binom Vr}$, where
\[
 X_e\coloneqq
 \begin{cases}
  1,&\text{if }e\text{ is an edge},\\
  0,&\text{otherwise}
 \end{cases}.
\]

We use the terminology of Austin~\cite[Section~3.9]{AustinExchangeable} and Dodos, Tyros, and Valettas~\cite[Section~1.2]{DTV}.

\begin{definition}\label{def:spreadable}
We call the random hypergraph \emph{spreadable} if its restrictions to finite vertex sets of the same size have the same distribution after increasing relabeling. More precisely, for any $a_1<\cdots<a_m$ and $b_1<\cdots<b_m$ in $V$ and every assignment of values $x_S\in\{0,1\}$ to the sets $S\in\binom{[m]}r$, we have
\[
 \Prb\left[X_{\{a_i\colon i\in S\}}=x_S\text{ for all }S\in\binom{[m]}r\right]
 =
 \Prb\left[X_{\{b_i\colon i\in S\}}=x_S\text{ for all }S\in\binom{[m]}r\right].
\]
\end{definition}
In particular, taking $m=r$ in \cref{def:spreadable} shows that every $r$-set has the same probability of being an edge; that is,
\[
 \Prb[X_e=1]=\Prb[X_{e'}=1]
 \quad\text{for all }e,e'\in\binom Vr.
\]

Let $(X_e)_{e\in\binom{\N}{r}}$ be a spreadable random $r$-graph with probability distribution $\nu$. For a finite ordered $r$-graph $F$, write $\Phi_\nu(F)$ for the probability that all edges of a fixed increasing copy of $F$ are present. More precisely, choose an increasing injection $\psi\colon V(F)\to\N$ and define
\[
 \Phi_\nu(F)\coloneqq\E\left[\prod_{e\in F}X_{\psi(e)}\right].
\]
\begin{fact}\label{fact:embedding-choice}
For every finite ordered $r$-graph $F$, the expectation defining $\Phi_\nu(F)$ does not depend on the choice of the increasing injection $\psi\colon V(F)\to\N$.
\end{fact}

\begin{proof}
Let $\psi,\psi'\colon V(F)\to\N$ be increasing injections. Since both preserve the order on $F$, spreadability gives the same joint distribution for the edge indicators on $\psi(V(F))$ and $\psi'(V(F))$, with vertices matched according to that order. Thus the probability that all edges of the corresponding copy of $F$ are present is the same for both injections.
\end{proof}
To formulate the independence assumption, define the \emph{$q$-shadow} of an $r$-graph $H$ by
\[
 \partial_qH\coloneqq\left\{S\in\binom{V(H)}q\colon S\subseteq e\text{ for some }e\in H\right\}.
\]

\begin{definition}\label{def:shadow-independent}
A random $r$-graph has \emph{$q$-shadow independence} if, whenever finite $r$-graphs $H_1,\ldots,H_t$ on subsets of its vertex set satisfy $\partial_qH_i\cap\partial_qH_{i'}=\varnothing$ for $i\ne i'$, the random vectors $(X_e)_{e\in H_i}$ are mutually independent.
\end{definition}

Note that the disjointness condition $\partial_qH_i\cap\partial_qH_{i'}=\varnothing$ is equivalent to $|e\cap e'|<q$ whenever $e\in H_i$, $e'\in H_{i'}$, and $i\ne i'$.

The following theorem shows that these invariance and independence assumptions force every ordered $q$-position-consistent hypergraph to occur with positive probability when the common edge probability is positive.

\begin{theorem}\label{lem:array-embedding}
Let $2\le q\le r-1$, and let $F$ be an ordered $q$-position-consistent $r$-graph with at least one edge and no isolated vertices. Write each edge as $e=\{v_1(e)<\cdots<v_r(e)\}$. For each $i\in[r]$, set $n_i\coloneqq\bigl|\{v_i(e)\colon e\in F\}\bigr|$. Let $(X_e)_{e\in\binom{\N}{r}}$ be a $q$-shadow-independent spreadable random $r$-graph with probability distribution $\nu$ and common edge probability $p$. Then
$\Phi_\nu(F)\ge p^{n_1\cdots n_r}$.
\end{theorem}

The next three lemmas establish the exchange and splitting properties and the multipartite bound needed for \cref{lem:array-embedding}.

For an ordered $r$-graph $G$ and distinct vertices $v,w\in V(G)$, write $G_{v\leftrightarrow w}$ for the ordered $r$-graph obtained by exchanging their positions in the vertex order and keeping the edge set unchanged.

\begin{lemma}\label{lem:adjacent-exchange}
Let $(X_e)_{e\in\binom{\N}{r}}$ be a spreadable random $r$-graph with probability distribution $\nu$. Suppose that $v<w$ are consecutive vertices of an ordered $r$-graph $G$ and no edge contains both. Then $\Phi_\nu(G)=\Phi_\nu(G_{v\leftrightarrow w})$.
\end{lemma}

\begin{proof}
For $M\ge1$, replace the consecutive pair $v,w$ in the vertex order of $G$ by a block of $2M+1$ new vertices. Number the resulting ordered set consecutively from $1$, and write $z_1<\cdots<z_{2M+1}$ for the integers assigned to the new block. Then every vertex preceding $v$ has image below $z_1$, and every vertex following $w$ has image above $z_{2M+1}$. Place $w$ at $z_{M+1}$, and identify all vertices other than $v$ with their assigned positions. For $t\in[2M+1]\setminus\{M+1\}$, let $A_t$ be the indicator that every edge containing $v$ is present when $v$ is placed at $z_t$, and let $Z$ be the indicator that all other edges are present. Thus
\[
 A_t\coloneqq\prod_{e\in G\colon v\in e}X_{(e\setminus\{v\})\cup\{z_t\}}
 \quad\text{and}\quad Z\coloneqq\prod_{e\in G\colon v\notin e}X_e.
\]

For $t\le M$, placing $v$ at $z_t$ preserves the order of $G$; for $t\ge M+2$, it gives the order of $G_{v\leftrightarrow w}$. In each case, $ZA_t$ is the indicator that all edges of the corresponding ordered copy are present. Define the left and right averages by $L_M\coloneqq M^{-1}\sum_{t=1}^M A_t$ and $R_M\coloneqq M^{-1}\sum_{t=M+2}^{2M+1}A_t$. By \cref{fact:embedding-choice},
\[
 \E[ZL_M]=\Phi_\nu(G)
 \quad\text{and}\quad
 \E[ZR_M]=\Phi_\nu(G_{v\leftrightarrow w}).
\]

Since $0\le Z\le1$, Cauchy--Schwarz gives
\[
 |\Phi_\nu(G)-\Phi_\nu(G_{v\leftrightarrow w})|
 =|\E[Z(L_M-R_M)]|\le\sqrt{\E[(L_M-R_M)^2]}.
\]
It therefore suffices to prove that $\E[(L_M-R_M)^2]\le2/M$.

Since no edge contains both $v$ and $w$, none of the edge indicators in $A_t$ involves $z_{M+1}$. The fixed vertices other than $w$ all lie outside the block. Thus the edge families occurring in $A_t$ are order-isomorphic for all candidate positions $t$. Likewise, the edge families occurring in $A_tA_{t'}$ are order-isomorphic for all distinct candidate positions $t,t'$, even when $z_t$ and $z_{t'}$ lie on opposite sides of $w$. By spreadability, $\E[A_t]=\mu$ and $\E[A_tA_{t'}]=c$ for all such $t,t'$ with $t\ne t'$, where $0\le c\le\mu\le1$.

Expanding $M^2L_M^2$ and $M^2R_M^2$ gives $M$ diagonal terms with expectation $\mu$ and $M(M-1)$ off-diagonal terms with expectation $c$ in each case. All $M^2$ terms in the expansion of $M^2L_MR_M$ have expectation $c$. After dividing by $M^2$, we obtain
\[
 \E[L_M^2]=\E[R_M^2]=\frac{\mu}{M}+\left(1-\frac1M\right)c
 \quad\text{and}\quad \E[L_MR_M]=c.
\]
Subtracting $2\E[L_MR_M]$ from $\E[L_M^2]+\E[R_M^2]$ gives
\begin{align*}
 \E[(L_M-R_M)^2]&=2\left[\frac{\mu}{M}+\left(1-\frac1M\right)c\right]-2c
 =\frac{2(\mu-c)}M\le\frac2M.
\end{align*}
Consequently, $|\Phi_\nu(G)-\Phi_\nu(G_{v\leftrightarrow w})|\le\sqrt{2/M}$ for every $M\ge1$. Letting $M\to\infty$ gives $\Phi_\nu(G)=\Phi_\nu(G_{v\leftrightarrow w})$.
\end{proof}

Let $G$ be an ordered $r$-graph, let $v\in V(G)$, and let $E_1,\ldots,E_\ell$ be a partition of $\{e\in G\colon v\in e\}$ into nonempty classes. For $1\le q\le r$, we call this partition \emph{$q$-shadow-disjoint} if $\partial_qE_s\cap\partial_qE_{s'}=\varnothing$ for all distinct $s,s'\in[\ell]$. Equivalently,
\begin{equation}\label{eq:split-condition}
 |e\cap e'|<q\quad\text{for all }e\in E_s\text{ and }e'\in E_{s'}\text{ with }s\ne s'.
\end{equation}

Fix the edge partition $E_1,\ldots,E_\ell$, and choose distinct new vertices $v_1,\ldots,v_\ell$ outside $V(G)$. For each permutation $\sigma$ of $[\ell]$, define $G^{\sigma,v}$ on the vertex set $(V(G)\setminus\{v\})\cup\{v_1,\ldots,v_\ell\}$ by
\[
 G^{\sigma,v}\coloneqq\{e\in G\colon v\notin e\}
 \cup\bigcup_{s\in[\ell]}\left\{(e\setminus\{v\})\cup\{v_s\}\colon e\in E_s\right\}.
\]
Its vertex order is obtained by replacing $v$ in the order of $G$ with the consecutive block $v_{\sigma(1)}<\cdots<v_{\sigma(\ell)}$, leaving the relative order of all other vertices unchanged. We call this operation \emph{splitting $v$ along the partition}.

\begin{lemma}\label{lem:splitting}
Let $2\le q\le r-1$, and let $(X_e)_{e\in\binom{\N}{r}}$ be a $q$-shadow-independent spreadable random $r$-graph with probability distribution $\nu$. Let $v$ be a vertex of an ordered $r$-graph $G$, and let $E_1,\ldots,E_\ell$ be a $q$-shadow-disjoint partition of the edges containing $v$. Then $\Phi_\nu(G)=\Phi_\nu(G^{\sigma,v})$ for every permutation $\sigma$ of $[\ell]$.
\end{lemma}

\begin{proof}
For $M\ge\max\{2,\ell\}$, replace $v$ in the vertex order of $G$ by a block of $M$ new vertices. Number the resulting ordered set consecutively from $1$, and write $z_1<\cdots<z_M$ for the integers assigned to the new block. Then every vertex preceding $v$ has image below $z_1$, and every vertex following $v$ has image above $z_M$. Identify the vertices of $V(G)\setminus\{v\}$ with their images. For $s\in[\ell]$ and $t\in[M]$, let $Y_{s,t}$ be the indicator that all edges in $E_s$ are present when $v$ is placed at $z_t$, and let $Z$ be the indicator that all edges not containing $v$ are present. Thus
\[
 Y_{s,t}\coloneqq\prod_{e\in E_s}X_{(e\setminus\{v\})\cup\{z_t\}}
 \quad\text{and}\quad Z\coloneqq\prod_{e\in G\colon v\notin e}X_e.
\]
We compare using the same candidate position for every edge class with choosing a candidate position independently for each class, allowing repetitions. Define
\[
 U_M\coloneqq\frac1M\sum_{t\in[M]}\prod_{s\in[\ell]}Y_{s,t}
 \quad\text{and}\quad V_M\coloneqq\prod_{s\in[\ell]}\left(\frac1M\sum_{t\in[M]}Y_{s,t}\right).
\]
Each summand in $U_M$ uses the same replacement vertex for all groups, so multiplying it by $Z$ gives the indicator that all edges of an ordered copy of $G$ are present.

No edge of a split hypergraph contains two copies of $v$, so \cref{lem:adjacent-exchange} makes all $\Phi_\nu(G^{\sigma,v})$ equal to a common value $b$. Expanding $V_M$ averages over independent uniform choices from $\{z_1,\ldots,z_M\}$ for the $\ell$ groups. Whenever the replacement vertices are distinct, the corresponding product, including $Z$, is the indicator that all edges of an ordered copy of some $G^{\sigma,v}$ are present and has expectation $b$ under $\nu$. The probability of a repeated choice is at most $\binom{\ell}{2}/M$. Since all products lie in $[0,1]$, it follows that
\[
 \E[ZU_M]=\Phi_\nu(G)
 \quad\text{and}\quad
 |\E[ZV_M]-b|\le\frac{\binom{\ell}{2}}M.
\]
The triangle inequality and Cauchy--Schwarz, with $0\le Z\le1$, therefore give
\[
 |\Phi_\nu(G)-b|
 \le\sqrt{\E[(U_M-V_M)^2]}+\frac{\binom{\ell}{2}}M.
\]
It remains to prove that $\E[(U_M-V_M)^2]\le1/M$.

For $s\ne s'$, $e\in E_s$, and $e'\in E_{s'}$, replacing $v$ by $z_t$ and $z_{t'}$, respectively, gives edges whose intersection has size $|e\cap e'|-1+\mathbf1\{t=t'\}\le|e\cap e'|<q$. Thus the edge families defining different rows have disjoint $q$-shadows, and \cref{def:shadow-independent} implies that the row vectors $(Y_{s,t})_{t=1}^M$, for $s\in[\ell]$, are mutually independent. All candidate positions have the same relative order with respect to the fixed vertices, so spreadability gives constants $\mu_s\coloneqq\E[Y_{s,t}]$ and $c_s\coloneqq\E[Y_{s,t}Y_{s,t'}]$ for $t\ne t'$, with $0\le c_s\le\mu_s\le1$.

For a fixed row $s$ and position $t$, averaging $\E[Y_{s,t}Y_{s,t'}]$ over $t'\in[M]$ gives $c_s+(\mu_s-c_s)/M$, since one choice has $t'=t$ and the other $M-1$ choices have $t'\ne t$. Averaging also over $t$ gives the same factor for the second moment of the row average. This factor also occurs in $\E[U_MV_M]$ because the common position in $U_M$ can be fixed before averaging the separate positions in $V_M$. Using the independence of the rows, and separating equal from distinct pairs of positions in $U_M^2$, we obtain
\begin{align*}
 \E[U_M^2] =\frac{1}{M}\prod_s\mu_s+(1-\frac{1}{M})\prod_s c_s\quad\text{and}\quad
 \E[V_M^2] =\E[U_MV_M]=\prod_s\left(c_s+\frac{\mu_s-c_s}{M}\right).
\end{align*}
Each factor in $\E[V_M^2]$ is at least $c_s$. Hence
\[
 \E[(U_M-V_M)^2]
 =\E[U_M^2]-\E[V_M^2]
 \le\frac1M\left(\prod_s\mu_s-\prod_s c_s\right)
 \le\frac1M.
\]
Substituting this into the comparison of occurrence probabilities gives
\[
 |\Phi_\nu(G)-b|\le \frac{1}{M^{1/2}}+\frac{\binom{\ell}{2}}M.
\]
Letting $M\to\infty$ gives $\Phi_\nu(G)=b=\Phi_\nu(G^{\sigma,v})$ for every permutation $\sigma$ of $[\ell]$.
\end{proof}

An ordered $r$-graph is \emph{position-partite} if its vertices are labeled by elements of $[r]$ so that the labels along every edge, in increasing vertex order, are $1,\ldots,r$. The vertices with each label form a \emph{position class}; these classes need not be intervals. The following bound uses only spreadability.

\begin{lemma}\label{lem:partite-count}
Let $(X_e)_{e\in\binom{\N}{r}}$ be a spreadable random $r$-graph with probability distribution $\nu$ and common edge probability $p$. If $G$ is position-partite with position classes of sizes $n_1,\ldots,n_r\ge1$, then
$\Phi_\nu(G)\ge p^{n_1\cdots n_r}$.
\end{lemma}

\begin{proof}
A consecutive pair of vertices whose labels are in decreasing order cannot lie in a common edge. By \cref{lem:adjacent-exchange}, we may exchange such pairs until the position classes occur as consecutive intervals in the order $1,\ldots,r$, preserving each edge's order and $\Phi_\nu(G)$. Adding every missing $r$-set containing one vertex from each position class gives a complete $r$-partite $r$-graph $K$ with $\Phi_\nu(G)\ge \Phi_\nu(K)$. Since $|K|=n_1\cdots n_r$, it suffices to show that $\Phi_\nu(K)\ge p^{|K|}$.

For $M\ge\max_i n_i$, choose intervals $B_1<\cdots<B_r$ of size $M$ and fix a realization of the random hypergraph. For positive integers $m_1,\ldots,m_r$, sample the vertices $x_{i,s}$, for $i\in[r]$ and $s\in[m_i]$, independently and uniformly from $B_i$, with replacement. Let $T_M(m_1,\ldots,m_r)$ be the probability, over these samples, that all crossing edges formed by one sampled vertex from each block are present. Thus
\[
 T_M(m_1,\ldots,m_r)\coloneqq\E_{\bx}\left[
 \prod_{(s_1,\ldots,s_r)\in[m_1]\times\cdots\times[m_r]}
 X_{\{x_{1,s_1},\ldots,x_{r,s_r}\}}\right].
\]
In particular, $D_M\coloneqq T_M(1,\ldots,1)$ is the crossing edge density on these blocks. The next step is to establish $T_M(n_1,\ldots,n_r)\ge D_M^{|K|}$.

Fix $i\in[r]$ and the samples outside $B_i$. Let $A_i$ denote the proportion of vertices $x\in B_i$ for which every edge formed by $x$ and one sampled vertex from each other block is present. Thus
\[
 A_i\coloneqq\frac1M\sum_{x\in B_i}
 \prod_{(s_h)_{h\ne i}\in\prod_{h\ne i}[m_h]}
 X_{\{x\}\cup\{x_{h,s_h}\colon h\ne i\}}.
\]
Each of the $m_i$ independent samples from $B_i$ must satisfy these edge conditions, so the conditional expectation of the product defining $T_M(m_1,\ldots,m_r)$ is $A_i^{m_i}$. Taking expectation over the samples outside $B_i$ and applying Jensen's inequality gives
\begin{align*}
 T_M(m_1,\ldots,m_r)
 =\E_{\bx}[A_i^{m_i}]
 \ge\bigl(\E_{\bx}[A_i]\bigr)^{m_i}
 =T_M(m_1,\ldots,m_{i-1},1,m_{i+1},\ldots,m_r)^{m_i}.
\end{align*}
Apply the preceding inequality successively for $i=1,\ldots,r$, reducing the $i$th argument of $T_M$ from $n_i$ to $1$ at step $i$. This gives $T_M(n_1,\ldots,n_r)\ge D_M^{|K|}$. This is the part-respecting form of Sidorenko's inequality for complete $r$-partite $r$-graphs; see Wigderson~\cite[Theorem~1 and the subsequent remark]{WigdersonSidorenko}. Taking expectation over the random hypergraph and using $\E[D_M]=p$, Jensen's inequality gives
\[
 \E[T_M(n_1,\ldots,n_r)]\ge\E[D_M^{|K|}]\ge(\E[D_M])^{|K|}=p^{|K|}.
\]
To compare $\E[T_M(n_1,\ldots,n_r)]$ with $\Phi_\nu(K)$, we bound the probability of repeated samples within a block. Whenever the sampled vertices are distinct within each block, sorting them within their blocks defines an increasing injection from $V(K)$ into $\N$. Spreadability implies that the expectation of the corresponding product over the random hypergraph is $\Phi_\nu(K)$. The probability of a repeated sample within some block is at most $M^{-1}\sum_i\binom{n_i}{2}$. Since all products lie in $[0,1]$,
\[
 \Phi_\nu(K)\ge\E[T_M(n_1,\ldots,n_r)]-\frac1M\sum_i\binom{n_i}{2}
 \ge p^{|K|}-\frac1M\sum_i\binom{n_i}{2}.
\]
Letting $M\to\infty$ gives $\Phi_\nu(K)\ge p^{|K|}$. Since $\Phi_\nu(G)\ge\Phi_\nu(K)$, this yields $\Phi_\nu(G)\ge p^{n_1\cdots n_r}$.
\end{proof}

\begin{proof}[Proof of \cref{lem:array-embedding}]
We use \cref{lem:splitting} to split each vertex of $F$ according to its positions in edges while preserving $\Phi_\nu(F)$. The resulting hypergraph is position-partite, so \cref{lem:partite-count} gives the required bound.

Write each edge of $F$ as $e=\{v_1(e)<\cdots<v_r(e)\}$. For a vertex $v$ and $i\in[r]$, put $E_i(v)\coloneqq\{e\in F\colon v_i(e)=v\}$. These groups satisfy
\[
 |e\cap e'|<q\quad\text{for all }e\in E_i(v)\text{ and }e'\in E_h(v)\text{ with }i\ne h.
\]
Indeed, otherwise $e\cap e'$ would contain a $q$-set containing $v$, and position consistency would force $v$ to occupy the same position in both edges.

Apply \cref{lem:splitting} at each vertex $v$ to the $q$-shadow-disjoint partition formed by the nonempty groups $E_i(v)$. Each split preserves the positions within every edge and can only decrease intersection sizes. Hence the groups at every vertex not yet split still satisfy \eqref{eq:split-condition}. The resulting $r$-graph $G$ is position-partite: its $i$th position class consists of the copies $(v,i)$ with $E_i(v)\ne\varnothing$, and has size $|\{v_i(e)\colon e\in F\}|=n_i$. By \cref{lem:splitting,lem:partite-count},
$\Phi_\nu(F)=\Phi_\nu(G)\ge p^{n_1\cdots n_r}$.
\end{proof}

For probability distributions $\mu$ and $\nu$ on the same finite set $\Omega$, their \emph{total variation distance} is defined by
\[
 \|\mu-\nu\|_{\mathrm{TV}}\coloneqq
 \frac{1}{2}\sum_{\omega\in\Omega}|\mu(\{\omega\})-\nu(\{\omega\})|
 =\max_{A\subseteq\Omega}|\mu(A)-\nu(A)|.
\]
Thus this distance is the largest difference between the probabilities assigned to the same event.

\emph{Remark.} The position-partite bound in \cref{lem:partite-count} can also be obtained from the distributional decomposition theorem of Dodos, Tyros, and Valettas~\cite[Theorem~1.4]{DTV}. Fix a position-partite $G$ with position classes of sizes $n_1,\ldots,n_r$, and fix $0<\delta<1$. Apply their theorem to the restriction of $\nu$ to a sufficiently large finite initial segment of $\N$, whose size is chosen in terms of $r$, $v(G)$, and $\delta$. The resulting finite mixture of kernel models with independent vertex variables approximates the distribution induced by $\nu$ on $[v(G)]$ within $\delta$ in total variation. Applying Jensen's inequality successively to the position classes within each model and then to the mixture weights, and letting $\delta\to0$ with $G$ fixed, gives $\Phi_\nu(G)\ge p^{n_1\cdots n_r}$. Together with \cref{lem:splitting}, this provides an alternative proof of \cref{lem:array-embedding}.

\section{Proof of the main theorem}\label{sec:zero}

We first introduce the partition notation and regularity tools needed to construct the finite random hypergraphs in \cref{lem:random-reduction}. The terminology follows R\"odl and Schacht~\cite[Section~2.3]{RSregular} and~\cite[Section~3.1.1]{RSremoval}.

Let $V_1,\ldots,V_t$ be disjoint vertex classes, called clusters, of equal size. For $S\subseteq[t]$, write $V_S\coloneqq\bigcup_{i\in S}V_i$. A set $e\subseteq V_{[t]}$ is \emph{crossing} if $|e\cap V_i|\le1$ for every $i\in[t]$. For partite hypergraphs, we consider only crossing cliques.

For an integer $k\ge2$, a \emph{family of partitions} through level $k-1$ consists of partitions of the crossing $\ell$-sets, for $1\le\ell\le k-1$, defined recursively as follows. For $2\le\ell<k$, crossing $\ell$-sets are first grouped according to the classes containing their $(\ell-1)$-subsets, and each group is then partitioned into labeled classes. The construction starts with the $1$-graphs $\binom{V_i}1$ for $i\in[t]$. Suppose its partition classes have been defined through level $\ell-1$, where $2\le\ell\le k$. Fix $S\in\binom{[t]}\ell$. For each $T\in\binom S{\ell-1}$, choose a partition class $P_T$ on $V_T$. Then put
\[
 B_S\coloneqq\bigcup_{T\in\binom S{\ell-1}}P_T.
\]
When $\mathcal K_\ell(B_S)\ne\varnothing$, the $(\ell-1)$-graph $B_S$ is called a \emph{polyad}. For $\ell<k$, partition each such $\mathcal K_\ell(B_S)$ into exactly $a_\ell$ nonempty classes, labeled $1,\ldots,a_\ell$~\cite[Definition~7]{RSregular}. These are the level-$\ell$ partition classes. Write $\ba\coloneqq(t,a_2,\ldots,a_{k-1})$ for the vector of partition parameters.

An \emph{address} $\alpha$ through level $h$ on $J\subseteq[t]$ records $J$ and the labels of selected classes $P_S$ on $V_S$, for $S\subseteq J$ with $2\le|S|\le h$~\cite[Definition~21]{RSremoval}. Here $1\le h\le\min\{k-1,|J|\}$, and we set $P_{\{i\}}\coloneqq\binom{V_i}1$ for $i\in J$. The classes are selected in increasing order of $|S|$. For each index set $S$, its label selects a class within the polyad determined by the classes already fixed on the $(|S|-1)$-subsets of $S$. Thus equal labels on different polyads need not refer to the same class. The selections must satisfy
\[
 P_S\subseteq\mathcal K_{|S|}\!\left(\bigcup_{T\in\binom S{|S|-1}}P_T\right)
 \quad\text{for every }S\subseteq J\text{ with }2\le|S|\le h.
\]
The crossing sets with address $\alpha$ form
\[
 \mathcal C(\alpha)\coloneqq\left\{e\in\binom{V_J}{|J|}\colon e\cap V_S\in P_S
 \text{ for every }S\subseteq J\text{ with }1\le|S|\le h\right\}.
\]
For $|J|=k$ and $h=k-1$, we have $\mathcal C(\alpha)=\mathcal K_k(B_\alpha)$, where $B_\alpha\coloneqq\bigcup_{S\in\binom J{k-1}}P_S$. $B_\alpha$ is a polyad whenever $\mathcal C(\alpha)\ne\varnothing$.

Let $B$ be an $\ell$-partite $(\ell-1)$-graph with $\mathcal K_\ell(B)\ne\varnothing$, and let $A\subseteq\mathcal K_\ell(B)$. For $0<\delta\le1$, $d\in[0,1]$, and an integer $R\ge1$, we call $A$ \emph{$(\delta,d,R)$-regular with respect to $B$}~\cite[Definition~11]{RSregular} if every union $U$ of at most $R$ sets $\mathcal K_\ell(B')$ with $B'\subseteq B$ and $|U|\ge\delta|\mathcal K_\ell(B)|$ satisfies
\[
 \left|\frac{|A\cap U|}{|U|}-d\right|\le\delta.
\]
We omit $R$ when $R=1$. Taking $U=\mathcal K_\ell(B)$ gives $|A|/|\mathcal K_\ell(B)|=d\pm\delta$.

An \emph{$h$-complex} on $u$ vertex classes is a collection $\mathcal H=(H^{(1)},\ldots,H^{(h)})$ of crossing hypergraphs, where $H^{(1)}=\bigcup_{i=1}^u\binom{V_i}1$ and $H^{(\ell)}\subseteq\mathcal K_\ell(H^{(\ell-1)})$ for $2\le\ell\le h$. For $\boldsymbol d=(d_2,\ldots,d_h)$, it is \emph{$(\eps,\boldsymbol d)$-regular} if $H^{(\ell)}[V_S]$ is $(\eps,d_\ell)$-regular with respect to $H^{(\ell-1)}[V_S]$ for every $\ell$ with $2\le\ell\le h$ and every $S\in\binom{[u]}\ell$~\cite[Definitions~1 and~5]{RSregular}.

R\"odl and Schacht~\cite[Theorem~2.3 and Lemma~4.1]{RSregular} proved a hypergraph regularity lemma allowing a prescribed initial partition. We use the following form, stated by Lin, Sun, Wang, and Zhou~\cite[Theorem~3.7]{LSWZ}.

\begin{theorem}[Regularity lemma]\label{thm:regularity}
Fix integers $s\ge k\ge2$, a constant $0<\sigma\le1$, and functions $\eps\colon\N^{k-1}\to(0,1]$ and $R\colon\N^{k-1}\to\N$. There are integers $T$ and $n_0$ such that the following holds. Let $H$ be a $k$-graph on $n\ge n_0$ vertices, where $sT!$ divides $n$, and fix an ordered equipartition of $V(H)$ into $s$ classes. Then $H$ admits a family of partitions through level $k-1$ with parameter vector $\ba=(t,a_2,\ldots,a_{k-1})$ such that
\begin{enumerate}
\item the vertex partition refines the initial partition, has $t$ classes of size $n/t$, and satisfies $s\le t\le T$; also, $a_\ell\le T$ for $2\le\ell\le k-1$;
\item every level-$\ell$ partition class is $(\eps(\ba),1/a_\ell)$-regular with respect to its $(\ell-1)$-polyad for $2\le\ell\le k-1$;
\item all but a $\sigma$ fraction of crossing $k$-sets lie in $\mathcal K_k(B)$ for a top polyad $B$ such that $H\cap\mathcal K_k(B)$ is $(\sigma,d_B,R(\ba))$-regular with respect to $B$, where $d_B=|H\cap\mathcal K_k(B)|/|\mathcal K_k(B)|$.
\end{enumerate}
\end{theorem}

The dense counting lemma of Kohayakawa, R\"odl, and Skokan~\cite[Theorem~6.5]{KRScounting} estimates the number of crossing cliques in a regular complex. We use the form stated by R\"odl and Schacht~\cite[Theorem~2.1]{RScounting}.

\begin{lemma}[Dense counting lemma]\label{lem:dense-counting}
For integers $2\le h\le u$ and constants $0<\mu<1$ and $d_0>0$, there are $\eps_0>0$ and $m_0$ with the following property. If $\mathcal H$ is an $(\eps_0,\boldsymbol d)$-regular $h$-complex on $u$ vertex classes of size $m\ge m_0$, where $d_\ell\ge d_0$ for $2\le\ell\le h$, then
\[
 |\mathcal K_u(H^{(h)})|=(1\pm\mu)m^u\prod_{\ell=2}^h d_\ell^{\binom u\ell}.
\]
\end{lemma}

To apply \cref{lem:dense-counting} to an address $\alpha$ through level $h$ on $J$, where $|J|=u$, take $H^{(\ell)}\coloneqq\bigcup_{S\in\binom J\ell}P_S$ for $1\le\ell\le h$ and $d_\ell\coloneqq1/a_\ell$ for $2\le\ell\le h$. Thus, for fixed $M$ and $\mu$, choosing $\eps(\ba)$ sufficiently small gives
\begin{equation}\label{eq:address-count}
 |\mathcal C(\alpha)|=(1\pm\mu)m^u\prod_{\ell=2}^{h}a_\ell^{-\binom u\ell}
\end{equation}
on $u\le M$ vertex classes of sufficiently large common size $m$. For $h=1$, we have $|\mathcal C(\alpha)|=m^u$. These formulas also justify choosing labels successively on any set of at most $M$ clusters. At a step corresponding to a subset $S$ with $2\le|S|\le k-1$, the classes already selected on its nonempty proper subsets form a regular complex. Applying \eqref{eq:address-count}, or the exact formula for $h=1$, shows that the corresponding polyad has a nonempty clique set. Each of its $a_{|S|}$ partition classes is therefore available for the next choice.

Cooley, Fountoulakis, K\"uhn, and Osthus~\cite[Theorem~2]{CFKO} proved an embedding lemma for regular hypergraph complexes. We use the following consequence, stated by Lin, Sun, Wang, and Zhou~\cite[Theorem~3.8]{LSWZ}, which allows the relative density to vary among the polyads corresponding to the edges of the hypergraph being embedded.

\begin{lemma}[Hypergraph embedding lemma]\label{lem:regular-embedding}
Let $k\ge2$, let $F$ be a $k$-graph on $[v(F)]$ with $v(F)\ge k$, and fix $d_*>0$. There is $\sigma_0>0$ with the following property. For every choice of $d_2,\ldots,d_{k-1}>0$ with $1/d_\ell\in\N$ for $2\le\ell\le k-1$, there are $\eps_0>0$ and integers $R_0$ and $m_0$ such that the following holds. The constants $\sigma_0$, $\eps_0$, $R_0$, and $m_0$ depend on $F$ only through $v(F)$.

Let $\mathcal H$ be an $(\eps_0,(d_2,\ldots,d_{k-1}))$-regular $(k-1)$-complex on vertex classes $V_1,\ldots,V_{v(F)}$ of size $m\ge m_0$. Let $H$ be a crossing $k$-graph on the same vertex classes. For each $e\in F$, put $B_e\coloneqq H^{(k-1)}[V_e]$ and suppose that $\mathcal K_k(B_e)\ne\varnothing$, that
\[
 d_e=\frac{|H\cap\mathcal K_k(B_e)|}{|\mathcal K_k(B_e)|}\ge d_*,
\]
and that $H\cap\mathcal K_k(B_e)$ is $(\sigma_0,d_e,R_0)$-regular with respect to $B_e$. Then $H$ contains a copy of $F$.
\end{lemma}

The regularity reduction will fix partition classes on $j$-subsets of the chosen clusters. To express the resulting restrictions on an ordered $r$-tuple, we use a density condition in which the relation imposed on the coordinates indexed by $I\in\binom{[r]}j$ may depend on $I$.

For an $r$-graph $H$ on $n$ vertices, let $\mathscr G=(G_I)_{I\in\binom{[r]}j}$ be a family of relations on ordered $j$-tuples of $V(H)$. Write $\mathcal K_r(\mathscr G)$ for the ordered $r$-tuples of distinct vertices whose restriction to each $I$ belongs to $G_I$. Let $e_H(\mathscr G)$ count those whose underlying $r$-set is an edge of $H$. The coordinatewise density condition requires that, for every such family $\mathscr G$,
\begin{equation}\label{eq:coordinate-density}
 e_H(\mathscr G)\ge d|\mathcal K_r(\mathscr G)|-\eta n^r.
\end{equation}
Reiher, R\"odl, and Schacht~\cite[Proposition~2.5]{RRS} proved that, for a fixed forbidden $r$-graph, this condition defines the same Tur\'an density as \eqref{eq:dense}. For $j\ge2$, one passes between the conditions by adjusting the error parameter $\eta$. For $j=1$, the argument first passes to an induced subhypergraph on a positive proportion of the vertices. Thus, if $\pi_{j,\unif}(F)>d$, there are arbitrarily large $F$-free $r$-graphs satisfying \eqref{eq:coordinate-density} with arbitrarily small $\eta$.

Fix $r,j,q$, and $N$ as in \cref{thm:zero}. The proof uses \cref{lem:array-embedding} for the zero-density direction and a random coloring for the converse and the lower bound.

\begin{proof}[Proof of \cref{thm:zero}]
Deleting isolated vertices preserves position consistency, and any copy of the remaining hypergraph extends to a copy of $F$ in every sufficiently large host. We may therefore assume that $F$ has no isolated vertices.

Suppose for a contradiction that $F$ has a $q$-position-consistent ordering but $\pi_{j,\unif}(F)>0$. Fix such an ordering and choose $0<d<\pi_{j,\unif}(F)$. We first use regularity to construct finite $F$-free random hypergraphs with $q$-shadow independence and edge probabilities bounded away from zero.

\begin{claim}\label{lem:random-reduction}
For every $M\ge\max\{r,v(F)\}$ there is a random $r$-graph $(X_e)_{e\in\binom{[M]}r}$ which is $F$-free almost surely, has $q$-shadow independence, and satisfies $\Prb[X_e=1]\ge d/2$ for every $e\in\binom{[M]}r$.
\end{claim}

\begin{proof}
We will construct the random $r$-graph on $M$ selected clusters of a regular partition of a large $F$-free host $H$. After fixing the labels through level $j$, we will sample the labels at levels $q,\ldots,r-1$ independently and retain an $r$-set when the corresponding $r$-partite subhypergraph of $H$ is regular with respect to its polyad and has relative density at least $d/8$. We will choose the clusters and their lower labels so that every $r$-set is retained with probability at least $d/2$.

\begin{enumerate}[ref=(\roman*)]
\item\label{item:reduction-setup}
Choose $0<\sigma\le\sigma_0$, where $\sigma_0$ is supplied by \cref{lem:regular-embedding} with $k=r$ and $d_*=d/8$, so that $16\sigma d^{-1}\binom Mr<\tfrac12$. This last bound will be used in part~\ref{item:reduction-selection} to choose clusters and lower labels so that, for every $r$-set of selected clusters, the probability of an irregular completion is at most $d/8$.

To apply the counting and embedding lemmas, choose $0<\mu\le1/100$ and functions $\eps(\ba)$ and $R(\ba)$ so that \cref{lem:dense-counting} applies on at most $M$ clusters and \cref{lem:regular-embedding} applies to $F$ with $d_*=d/8$ and $d_\ell=1/a_\ell$ for $2\le\ell\le r-1$. Let $T$ be the upper bound on $t,a_2,\ldots,a_{r-1}$ supplied by \cref{thm:regularity} for an initial partition into $M$ classes. There are only finitely many possible parameter vectors $\ba$ with entries at most $T$. We will choose the number $n$ of vertices of the host so that $n/T$ exceeds the cluster-size thresholds in both lemmas for every such vector, ensuring that all these applications are valid for the family of partitions supplied by \cref{thm:regularity}.

For every sufficiently large $n$ and any family of partitions supplied by \cref{thm:regularity} on a host with $n$ vertices, we have $m=n/t$ and $t,a_2,\ldots,a_j\le T$. Thus \eqref{eq:address-count} gives $|\mathcal C(\alpha)|\ge c_Tn^r$ for every address $\alpha$ through level $j$ on $r$ clusters, with a uniform constant $c_T>0$. Once $T$ is known, choose $0<\eta\le dc_T/8$.

Using the equality of the Tur\'an densities associated with \eqref{eq:dense} and \eqref{eq:coordinate-density}, take an $F$-free $r$-graph satisfying \eqref{eq:coordinate-density} with error $\eta/2$ and of sufficiently large order to meet the preceding size requirements after removing fewer than $MT!$ vertices. Delete fewer than $MT!$ vertices to obtain an $r$-graph $H$ whose order $n$ is divisible by $MT!$.

The original host has fewer than $n+MT!$ vertices, and $(\eta/2)(n+MT!)^r\le\eta n^r$ for sufficiently large $n$. Extend each test relation on $V(H)$ to the original vertex set by allowing no tuples involving deleted vertices. Applying the original density inequality to these extensions gives \eqref{eq:coordinate-density} for $H$ with error $\eta$.

Apply \cref{thm:regularity} to $H$ with $k=r$, $s=M$, and an initial equipartition into $M$ classes to obtain a family of partitions through level $r-1$ with parameter vector $\ba=(t,a_2,\ldots,a_{r-1})$, where $M\le t\le T$ and the clusters have common size $m=n/t$.

\item\label{item:reduction-model}
For the partition obtained in part~\ref{item:reduction-setup}, fix any $M$ clusters, order them increasingly, and identify them with $[M]$. Fix also any choice of labels on their subsets of sizes $2,\ldots,j$.

With these clusters and lower labels fixed, sample independent uniform labels $\xi_S\in[a_{|S|}]$ for every $S\subseteq[M]$ with $q\le|S|\le r-1$. Starting with the fixed classes through level $j$, process these subsets in increasing order of size. For each $S$, let $P_S$ be the partition class with label $\xi_S$ within the polyad determined by the classes already selected on the $(|S|-1)$-subsets of $S$.

For $e\in\binom{[M]}r$, its top polyad is $B_e\coloneqq\bigcup_{S\in\binom e{r-1}}P_S$; put $d_e\coloneqq|H\cap\mathcal K_r(B_e)|/|\mathcal K_r(B_e)|$. Declare $e$ to be an edge of the random $r$-graph precisely when $H\cap\mathcal K_r(B_e)$ is regular with respect to $B_e$ and has density at least $d/8$. Equivalently, define
\[
 X_e\coloneqq\mathbf1\left\{H\cap\mathcal K_r(B_e)\text{ is }(\sigma,d_e,R(\ba))\text{-regular with respect to }B_e\text{ and }d_e\ge\frac{d}{8}\right\}.
\]

\item\label{item:reduction-properties}
We first verify $q$-shadow independence and $F$-freeness for every fixed choice of clusters and lower labels. The only random inputs to $X_e$ are the labels $\xi_S$ with $S\subseteq e$ and $q\le|S|\le r-1$. Edge families with disjoint $q$-shadows use disjoint sets of these labels. The vectors of edge indicators indexed by these families are therefore mutually independent, as required by \cref{def:shadow-independent}.

Suppose a realization admits an injective map $\phi\colon V(F)\to[M]$ with $X_{\phi(e)}=1$ for every $e\in F$. Assign each vertex $v$ of $F$ to the selected cluster indexed by $\phi(v)$. Together with the complete $1$-graphs on the corresponding clusters, the partition classes selected on subsets of $\phi(V(F))$ of sizes $2,\ldots,r-1$ form an $(r-1)$-complex. This complex is $(\eps(\ba),(1/a_2,\ldots,1/a_{r-1}))$-regular. For every edge $e\in F$, the $r$-graph $H\cap\mathcal K_r(B_{\phi(e)})$ is $(\sigma,d_{\phi(e)},R(\ba))$-regular with respect to $B_{\phi(e)}$ and has relative density $d_{\phi(e)}\ge d/8$.

By the choices of $\eps(\ba)$, $R(\ba)$, and the number of vertices of the host, we have $\eps(\ba)\le\eps_0$, $\sigma\le\sigma_0$, $R(\ba)\ge R_0$, and $m\ge m_0$ for the parameters in \cref{lem:regular-embedding}. Thus that lemma embeds $F$ in $H$, a contradiction. Every realization is therefore $F$-free.

\item\label{item:reduction-density}
It remains to choose the clusters and their labels through level $j$ so that $\Prb[X_e=1]\ge d/2$ for every $e\in\binom{[M]}r$. In this part we establish an average-density bound valid for every lower address; part~\ref{item:reduction-selection} will select clusters and lower labels so that the probability of an irregular completion is small for every $r$-set of selected clusters.

Fix any $r$ clusters of the partition and any address $\alpha$ through level $j$ on them. For $j\ge2$ and $I\in\binom{[r]}j$, let $G_I$ consist of the tuples obtained by listing, in cluster order, the vertices of each edge in the partition class selected at positions $I$; for $j=1$, let $G_{\{i\}}$ be the cluster in position $i$. Listing vertices in cluster order gives a bijection between $\mathcal C(\alpha)$ and $\mathcal K_r(\mathscr G)$. Thus \eqref{eq:coordinate-density} gives
\[
 \frac{|H\cap\mathcal C(\alpha)|}{|\mathcal C(\alpha)|}
 \ge d-\frac{\eta n^r}{|\mathcal C(\alpha)|}\ge\frac{7d}{8}.
\]
Extend $\alpha$ by choosing labels independently and uniformly on all subsets of sizes $q,\ldots,r-1$ of the $r$ clusters specified by $\alpha$. Let $\omega$ denote the completed address, and put $d_\omega\coloneqq|H\cap\mathcal C(\omega)|/|\mathcal C(\omega)|$, whether or not $H\cap\mathcal C(\omega)$ is regular with respect to the corresponding polyad. To bound $\E[d_\omega\mid\alpha]$, first note that the sets $\mathcal C(\omega)$ partition $\mathcal C(\alpha)$, so
\[
 \frac{|H\cap\mathcal C(\alpha)|}{|\mathcal C(\alpha)|}
 =\frac{\sum_\omega|\mathcal C(\omega)|d_\omega}
        {\sum_\omega|\mathcal C(\omega)|},
\]
where both sums range over all completions $\omega$ of $\alpha$. This identity weights each density $d_\omega$ by $|\mathcal C(\omega)|$, whereas the random choice of labels gives equal probability to each completion. By \eqref{eq:address-count}, the values $|\mathcal C(\omega)|$, as $\omega$ ranges over the completions of $\alpha$, differ by a factor of at most $(1+\mu)/(1-\mu)$, so
\begin{equation}\label{eq:mean-top-density}
 \E[d_\omega\mid\alpha]\ge\frac{1-\mu}{1+\mu}\frac{7d}{8}\ge\frac{3d}{4}.
\end{equation}
Since $\alpha$ was arbitrary, \eqref{eq:mean-top-density} gives $\E[d_e]\ge3d/4$ for every candidate edge in every model from part~\ref{item:reduction-model}.

\item\label{item:reduction-selection}
To find such a choice of clusters and lower labels, we first bound the proportion of irregular completed addresses over the whole partition. Let $\Omega$ be the set of addresses through level $r-1$ on all $r$-subsets of the $t$ clusters, and let $\Omega_{\mathrm{irr}}\subseteq\Omega$ consist of those addresses $\omega$ for which $H\cap\mathcal C(\omega)$ is not $(\sigma,d_\omega,R(\ba))$-regular with respect to the corresponding top polyad. For $\omega\in\Omega$, write $w_\omega\coloneqq|\mathcal C(\omega)|$. The regularity lemma gives
\[
 \sum_{\omega\in\Omega_{\mathrm{irr}}}w_\omega
 \le \sigma\sum_{\omega\in\Omega}w_\omega.
\]
Put $w_{\min}\coloneqq\min_{\omega\in\Omega}w_\omega$ and $w_{\max}\coloneqq\max_{\omega\in\Omega}w_\omega$. By \eqref{eq:address-count}, $\frac{w_{\max}}{w_{\min}}\le\frac{1+\mu}{1-\mu}<2$. 
Consequently,
\[
 \frac{|\Omega_{\mathrm{irr}}|}{|\Omega|}
 \le \frac{w_{\max}}{w_{\min}}
      \frac{\sum_{\omega\in\Omega_{\mathrm{irr}}}w_\omega}
           {\sum_{\omega\in\Omega}w_\omega}
 \le 2\sigma.
\]
For this selection argument, choose $M$ of the $t$ clusters uniformly and order them increasingly, then choose labels independently and uniformly on all subsets of sizes $2,\ldots,j$ of the selected $M$ clusters. For $e\in\binom{[M]}r$, let $b_e$ be the conditional probability that completing the address fixed on the clusters indexed by $e$ gives an address in $\Omega_{\mathrm{irr}}$; that is,
\[
 b_e\coloneqq\Prb\left[
  \omega\in\Omega_{\mathrm{irr}}
  \,\middle|\,
  \text{the selected clusters and labels through level }j
 \right],
\]
where $\omega$ is obtained by completing the address on the clusters indexed by $e$ with independent uniform labels at levels $q,\ldots,r-1$.

Every $r$-subset of the $t$ clusters is contained in the selected set of $M$ clusters with probability $\binom Mr/\binom tr$. Each $r$-subset supports exactly $\prod_{\ell=2}^{r-1}a_\ell^{\binom r\ell}$ completed addresses, since every $\ell$-subset of cluster indices has $a_\ell$ available labels. For each selected $r$-subset of clusters, the completed address is uniformly distributed over these addresses. Using this inclusion probability and the fact that every $r$-subset of clusters supports the same number of completed addresses, we obtain, after averaging over the selected clusters and their labels through level $j$,
\[
 \E\left[\sum_{e\in\binom{[M]}r}b_e\right]
 =\binom Mr\frac{|\Omega_{\mathrm{irr}}|}{|\Omega|}
 \le2\sigma\binom Mr.
\]
Since $\max_e b_e\le\sum_e b_e$, Markov's inequality gives
\[
 \Prb\left[\max_e b_e>\frac{d}{8}\right]\le\frac{8}{d}\E\left[\sum_e b_e\right]
 \le\frac{16\sigma}{d}\binom Mr<\frac{1}{2}.
\]

Fix a choice of clusters and labels through level $j$ with $b_e\le d/8$ for every $e$. Use this choice in the construction of part~\ref{item:reduction-model}, with only the labels at levels $q,\ldots,r-1$ random. Part~\ref{item:reduction-properties} gives $q$-shadow independence and $F$-freeness. The bound \eqref{eq:mean-top-density} from part~\ref{item:reduction-density} holds for every lower address, so $\E[d_e]\ge3d/4$ for every $e$ under this law.

The outcomes with $X_e=0$ contribute at most $b_e+d/8$ to $\E[d_e]$. Since $d_e\le1$ on the event $X_e=1$, we obtain
\[
 \Prb[X_e=1]\ge\E[d_e]-b_e-\frac{d}{8}
 \ge\frac{3d}{4}-\frac{d}{8}-\frac{d}{8}=\frac{d}{2}.\qedhere
\]
\end{enumerate}
\end{proof}

To apply \cref{lem:array-embedding}, we next pass from these finite random hypergraphs to a spreadable limit.

\begin{claim}\label{clm:spreadable-limit}
There is a $q$-shadow-independent spreadable random $r$-graph on $\N$ that is $F$-free almost surely and has common edge probability $p\ge d/2$.
\end{claim}

\begin{proof}
We first use Ramsey's theorem to construct finite random $r$-graphs whose laws are approximately invariant under increasing relabelings. We then take a diagonal limit and verify its properties in parts~\ref{item:limit-spreadability} and~\ref{item:limit-preservation}.

\begin{enumerate}[ref=(\roman*)]
\item\label{item:limit-finite-models}
Fix an integer $n\ge r$. For a finite random $r$-graph on an ordered vertex set and a $k$-set $S=\{s_1<\cdots<s_k\}$ of its vertices, let $\nu_S$ be the distribution of $\bigl(X_{\{s_i\colon i\in I\}}\bigr)_{I\in\binom{[k]}r}$. Thus $\nu_S$ is the induced distribution after relabeling $S$ increasingly as $[k]$.

We seek $n$ vertices such that, for each $r\le k\le n$, the induced distributions on any two $k$-sets differ in total variation by at most $1/n$. For each $r\le k\le n$, partition the compact simplex of probability distributions on $\{0,1\}^{\binom{[k]}r}$ into $c_k$ parts of total variation diameter at most $1/n$. Color a $k$-set $S$ by the part containing $\nu_S$. The numbers $c_k$ are fixed before the finite random hypergraph is chosen.

Write $R_k(m;c)$ for the least integer $L$ such that every $c$-coloring of $\binom{[L]}k$ has a monochromatic $m$-set. Ramsey's theorem~\cite{Ramsey} ensures that these numbers are finite. Set $m_{n+1}\coloneqq n$ and define, in reverse order,
\[
 m_k\coloneqq R_k(m_{k+1};c_k)
 \quad\text{for }k=n,n-1,\ldots,r.
\]
Choose $M\ge\max\{r,v(F),m_r\}$ and apply \cref{lem:random-reduction} on $[M]$.

In this finite model, apply Ramsey's theorem successively to the colorings of the $k$-sets for $k=r,r+1,\ldots,n$. At stage $k$, at least $m_k=R_k(m_{k+1};c_k)$ vertices are available, so retain a monochromatic subset of size $m_{k+1}$. Denote the final set by $A_n$, so that $|A_n|=m_{n+1}=n$. Each restriction preserves monochromaticity at every earlier value of $k$. Thus, for each $r\le k\le n$, all $k$-sets in $A_n$ have their distributions in the same part, so
\begin{equation}\label{eq:ramsey-approximation}
 \|\nu_S-\nu_T\|_{\mathrm{TV}}\le\frac1n
 \quad\text{for }r\le k\le n\text{ and all }S,T\in\binom{A_n}k.
\end{equation}

The coloring and the choice of $A_n$ depend only on the probability law of the finite model. Relabel $A_n$ increasingly as $[n]$, and let $\nu_n$ be the distribution of the induced random $r$-graph. Restriction and relabeling preserve $F$-freeness, $q$-shadow independence, and the lower bound $d/2$ on each edge probability.

\item\label{item:limit-law}
Carry out part~\ref{item:limit-finite-models} for every $n\ge r$. Fix $k\ge r$. As $n\to\infty$ through $n\ge k$, the restrictions of $\nu_n$ to $[k]$ lie in a compact finite-dimensional simplex. Successively select subsequences converging in total variation for $k=r,r+1,\ldots$, and take a diagonal subsequence along which all these restrictions converge. In the rest of this argument, $n$ ranges over that subsequence.

For each fixed $k$, every probability in the marginal distribution on $[k]$ is a finite sum of probabilities of edge-indicator configurations on $[k+1]$. Passing to the limit in these identities shows that the limiting distribution on $[k]$ is the marginal of the limiting distribution on $[k+1]$. These consistent laws define a probability distribution $\nu$ on $\{0,1\}^{\binom{\N}r}$.

\item\label{item:limit-spreadability}
To verify spreadability of the law $\nu$ from part~\ref{item:limit-law}, fix a $k$-set $S\subseteq\N$ with $k\ge r$. For all sufficiently large $n$, the set $S$ is contained in $[n]$. Under $\nu_n$, the distributions induced on $S$ and $[k]$, after increasing relabeling onto $[k]$, have total variation distance at most $1/n$ by \eqref{eq:ramsey-approximation}. Both laws converge, since $S\cup[k]$ lies in a fixed initial segment of $\N$, and their limits are therefore equal. Invariance for sets of fewer than $r$ vertices is automatic, since they contain no $r$-sets. Hence the random $r$-graph with distribution $\nu$ is spreadable.

For each fixed $e\in\binom{\N}r$, the bound $\Prb_{\nu_n}[X_e=1]\ge d/2$ holds for all sufficiently large $n$ and passes to the limit. Spreadability makes these probabilities equal to a common value $p\ge d/2$.

\item\label{item:limit-preservation}
With spreadability and the edge-probability bound established in part~\ref{item:limit-spreadability}, it remains to verify $F$-freeness and $q$-shadow independence.

For any fixed copy of $F$ on $\N$, the event that all its edges are present depends on finitely many edge indicators. Its probability is zero under every sufficiently large finite model, so finite-dimensional convergence makes its probability under $\nu$ zero as well. There are only countably many such copies on $\N$. The probability that any of them is present is therefore zero, and the limiting random hypergraph is $F$-free almost surely.

Finally, fix finite $r$-graphs $H_1,\ldots,H_t$ on subsets of $\N$ with pairwise disjoint $q$-shadows and prescribed values $x_e\in\{0,1\}$ for $e\in H_1\cup\cdots\cup H_t$. For all sufficiently large $n$, the finite models satisfy
\[
 \Prb_{\nu_n}[X_e=x_e\text{ for all }e\in H_1\cup\cdots\cup H_t]
 =\prod_{i=1}^t\Prb_{\nu_n}[X_e=x_e\text{ for all }e\in H_i].
\]
All edges indexing the indicators in this identity are contained in a fixed finite vertex set. Each probability in the identity therefore converges to the corresponding probability under $\nu$. Taking limits in the finite product gives the same factorization with $\nu_n$ replaced by $\nu$. Since the $H_i$ and the values $x_e$ were arbitrary, the limiting random $r$-graph has $q$-shadow independence.\qedhere
\end{enumerate}
\end{proof}

Let $\nu$ be the probability distribution supplied by \cref{clm:spreadable-limit}, with common edge probability $p\ge d/2$. Then \cref{lem:array-embedding}, with $n_i$ as defined there, gives the contradiction
\[
 0=\Phi_\nu(F)\ge p^{n_1\cdots n_r}>0.
\]
This proves the zero-density direction.

For the converse and the bound $N^{-N}$, we use the following construction, extending the random coloring construction of Reiher, R\"odl, and Schacht~\cite[Fact~1.3]{RRSzero} to $q$-sets.

\begin{claim}\label{clm:coloring-construction}
For every $\eta>0$ and all sufficiently large $n$, there is an $(N^{-N},\eta,j)$-dense $r$-graph on $[n]$ whose natural ordering is $q$-position-consistent.
\end{claim}

\begin{proof}
The color of a $q$-set specifies the positions it must occupy in any retained edge containing it. Assign each $q$-set of $[n]$ an independent uniform color in $\binom{[r]}q$, and denote the coloring by $\chi$. For $e=\{v_1<\cdots<v_r\}$, write $e_I\coloneqq\{v_i\colon i\in I\}$ for $I\in\binom{[r]}q$ and put
\[
 H_\chi\coloneqq\left\{e\in\binom{[n]}r\colon
 \chi(e_I)=I\text{ for every }I\in\binom{[r]}q\right\}.
\]
If $e_I=e'_J$ for edges $e$ and $e'$ in $H_\chi$ and $I,J\in\binom{[r]}q$, then $I=\chi(e_I)=\chi(e'_J)=J$. Hence the natural order of $H_\chi$ is $q$-position-consistent.

Fix $\eta>0$. We verify the density condition simultaneously for all $j$-graphs $G$ on $[n]$. For a fixed $j$-graph $G$ on $[n]$, put $Z_G\coloneqq|H_\chi\cap\mathcal K_r(G)|$. Each $r$-set imposes $N$ independent color requirements, each of probability $1/N$. Thus
\[
 \Prb[e\in H_\chi]=N^{-N}
 \quad\text{and}\quad \E[Z_G]=N^{-N}|\mathcal K_r(G)|.
\]
Recoloring a fixed $q$-set $S$ can change whether $e\in H_\chi$ only when $S\subseteq e$, so it changes $Z_G$ by at most $\binom{n-q}{r-q}\le n^{r-q}$. The bounded-differences inequality~\cite{McDiarmid} gives
\[
 \Prb\bigl[|Z_G-\E[Z_G]|>\eta n^r\bigr]
 \le2\exp\left(-\frac{2\eta^2n^{2r}}{\binom nq n^{2r-2q}}\right)
 \le2\exp(-2\eta^2n^q).
\]
To obtain the density bound simultaneously for all $j$-graphs $G$, apply the union bound over the $2^{\binom nj}$ choices. Since $q=j+1$, this gives
\[
 \Prb\bigl[\exists G\colon |Z_G-\E[Z_G]|>\eta n^r\bigr]
 \le2^{1+\binom nj}\exp(-2\eta^2n^{j+1})=o(1)
 \quad\text{as }n\to\infty.
\]
For all sufficiently large $n$, some realization of $H_\chi$ is therefore $(N^{-N},\eta,j)$-dense and has a $q$-position-consistent natural ordering.
\end{proof}

Any copy of $F$ in a hypergraph supplied by \cref{clm:coloring-construction} would induce a $q$-position-consistent ordering of $F$ by pulling back the natural order on $[n]$ along the embedding. If $F$ has no such ordering, these hypergraphs are therefore $F$-free, giving $\pi_{j,\unif}(F)\ge N^{-N}$. This proves the converse and the lower bound.

It remains to deduce the conclusion for forbidden families. If $\F$ contains a member $F$ with a $q$-position-consistent ordering, then
\[
 0\le\pi_{j,\unif}(\F)\le\pi_{j,\unif}(F)=0.
\]
Otherwise, every hypergraph supplied by \cref{clm:coloring-construction} avoids all members of $\F$ simultaneously, even if $\F$ is infinite. Hence $\pi_{j,\unif}(\F)\ge N^{-N}$.
\end{proof}

\section*{Declaration on the use of AI}
The authors used generative AI tools to assist in discussing proof strategies, checking proofs, and improving exposition.

\begingroup
\setlength{\parskip}{3pt plus 1pt}
\bibliographystyle{abbrv}
\bibliography{zero_uniform_density}

\begin{thebibliography}{10}

\bibitem{AustinExchangeable}
T.~Austin.
\newblock On exchangeable random variables and the statistics of large graphs
  and hypergraphs.
\newblock {\em Probab. Surv.}, 5:80--145, 2008.
\newblock \href{https://doi.org/10.1214/08-PS124}{doi:10.1214/08-PS124}.

\bibitem{BPRRSQuasirandomLinks}
S.~Berger, S.~Piga, C.~Reiher, V.~R{\"o}dl, and M.~Schacht.
\newblock {Tur{\'a}n} density of cliques of order five in {$3$}-uniform
  hypergraphs with quasirandom links.
\newblock {\em Trans. Amer. Math. Soc.}, 378(12):8283--8318, 2025.
\newblock \href{https://doi.org/10.1090/tran/9174}{doi:10.1090/tran/9174}.

\bibitem{BucicTetrahedron}
M.~Buci{\'c}.
\newblock The uniform {Tur{\'a}n} density of the tetrahedron.
\newblock Preprint, 2026.
\newblock \href{https://arxiv.org/abs/2609.11802v2}{arXiv:2609.11802v2}.

\bibitem{CFKO}
O.~Cooley, N.~Fountoulakis, D.~K{\"u}hn, and D.~Osthus.
\newblock Embeddings and {Ramsey} numbers of sparse {$k$}-uniform hypergraphs.
\newblock {\em Combinatorica}, 29(3):263--297, 2009.
\newblock
  \href{https://doi.org/10.1007/s00493-009-2356-y}{doi:10.1007/s00493-009-2356-y}.

\bibitem{DLLWYCodegree}
L.~Ding, A.~Lamaison, H.~Liu, S.~Wang, and H.~Yang.
\newblock On {$3$}-graphs with vanishing codegree {Tur{\'a}n} density.
\newblock {\em J. Lond. Math. Soc. (2)}, 112(3):e70281, 2025.
\newblock \href{https://doi.org/10.1112/jlms.70281}{doi:10.1112/jlms.70281}.

\bibitem{DTV}
P.~Dodos, K.~Tyros, and P.~Valettas.
\newblock Decompositions of finite high-dimensional random arrays.
\newblock {\em Fund. Math.}, 268(2):101--150, 2025.
\newblock
  \href{https://doi.org/10.4064/fm221004-6-11}{doi:10.4064/fm221004-6-11}.

\bibitem{ErdosExtremal}
P.~Erd{\H{o}}s.
\newblock On extremal problems of graphs and generalized graphs.
\newblock {\em Israel J. Math.}, 2(3):183--190, 1964.
\newblock \href{https://www.renyi.hu/~p_erdos/1964-13.pdf}{Author reprint}.

\bibitem{ErdosSimonovits}
P.~Erd{\H{o}}s and M.~Simonovits.
\newblock A limit theorem in graph theory.
\newblock {\em Studia Sci. Math. Hungar.}, 1:51--57, 1966.
\newblock \href{https://www.renyi.hu/~miki/ErdosSimLim.pdf}{Author reprint}.

\bibitem{ErdosSos}
P.~Erd{\H{o}}s and V.~T. S{\'o}s.
\newblock On {Ramsey--Tur{\'a}n} type theorems for hypergraphs.
\newblock {\em Combinatorica}, 2(3):289--295, 1982.
\newblock \href{https://doi.org/10.1007/BF02579235}{doi:10.1007/BF02579235}.

\bibitem{ErdosStone}
P.~Erd{\H{o}}s and A.~H. Stone.
\newblock On the structure of linear graphs.
\newblock {\em Bull. Amer. Math. Soc.}, 52(12):1087--1091, 1946.
\newblock
  \href{https://doi.org/10.1090/S0002-9904-1946-08715-7}{doi:10.1090/S0002-9904-1946-08715-7}.

\bibitem{GKLMinimum}
F.~Garbe, D.~Kr{\'a}{\v{l}}, and A.~Lamaison.
\newblock Hypergraphs with minimum positive uniform {Tur{\'a}n} density.
\newblock {\em Israel J. Math.}, 259(2):701--726, 2024.
\newblock
  \href{https://doi.org/10.1007/s11856-023-2554-0}{doi:10.1007/s11856-023-2554-0}.

\bibitem{GKV}
R.~Glebov, D.~Kr{\'a}{\v{l}}, and J.~Volec.
\newblock A problem of {Erd{\H{o}}s and S{\'o}s} on {$3$}-graphs.
\newblock {\em Israel J. Math.}, 211(1):349--366, 2016.
\newblock
  \href{https://doi.org/10.1007/s11856-015-1267-4}{doi:10.1007/s11856-015-1267-4}.

\bibitem{KeevashSurvey}
P.~Keevash.
\newblock Hypergraph {Tur{\'a}n} problems.
\newblock In R.~Chapman, editor, {\em Surveys in Combinatorics, 2011}, volume
  392 of {\em London Math. Soc. Lecture Note Ser.}, pages 83--140. Cambridge
  Univ. Press, Cambridge, 2011.
\newblock
  \href{https://doi.org/10.1017/CBO9781139004114.004}{doi:10.1017/CBO9781139004114.004}.

\bibitem{KKLLSWTetrahedron}
B.~Kielak, D.~Kr{\'a}{\v{l}}, A.~Lamaison, H.~Liu, X.~Shu, and Z.~Wu.
\newblock Solution of uniform {Tur{\'a}n}'s {Tetrahedron Problem}.
\newblock Preprint, 2026.
\newblock \href{https://arxiv.org/abs/2609.08336v3}{arXiv:2609.08336v3}.

\bibitem{KNRSlinear}
Y.~Kohayakawa, B.~Nagle, V.~R{\"o}dl, and M.~Schacht.
\newblock Weak hypergraph regularity and linear hypergraphs.
\newblock {\em J. Combin. Theory Ser. B}, 100(2):151--160, 2010.
\newblock
  \href{https://doi.org/10.1016/j.jctb.2009.05.005}{doi:10.1016/j.jctb.2009.05.005}.

\bibitem{KRScounting}
Y.~Kohayakawa, V.~R{\"o}dl, and J.~Skokan.
\newblock Hypergraphs, quasi-randomness, and conditions for regularity.
\newblock {\em J. Combin. Theory Ser. A}, 97(2):307--352, 2002.
\newblock
  \href{https://doi.org/10.1006/jcta.2001.3217}{doi:10.1006/jcta.2001.3217}.

\bibitem{LamaisonPalettes}
A.~Lamaison.
\newblock Palettes determine uniform {Tur{\'a}n} density.
\newblock Preprint, 2024.
\newblock \href{https://arxiv.org/abs/2408.09643}{arXiv:2408.09643}.

\bibitem{LamaisonBeyond}
A.~Lamaison.
\newblock Uniform {Tur{\'a}n} density beyond {$3$}-graphs.
\newblock Preprint, 2025.
\newblock \href{https://arxiv.org/abs/2508.20696v1}{arXiv:2508.20696v1}.

\bibitem{LMZVanishing}
H.~Lin, J.~Ma, and W.~Zhou.
\newblock Vanishing uniform {Tur{\'a}n} densities: Characterization and jumps.
\newblock Preprint, 2026.
\newblock \href{https://arxiv.org/abs/2610.28574v1}{arXiv:2610.28574v1}.

\bibitem{LSWZ}
H.~Lin, G.~Sun, G.~Wang, and W.~Zhou.
\newblock Uniform {Tur{\'a}n} densities of {$k$}-uniform hypergraphs.
\newblock Preprint, 2026.
\newblock \href{https://arxiv.org/abs/2605.15105v1}{arXiv:2605.15105v1}.

\bibitem{LWZMinimum}
H.~Lin, G.~Wang, and W.~Zhou.
\newblock The minimum positive uniform {Tur{\'a}n} density in uniformly dense
  {$k$}-uniform hypergraphs.
\newblock {\em SIAM J. Discrete Math.}, 39(4):2350--2378, 2025.
\newblock \href{https://doi.org/10.1137/24M1658218}{doi:10.1137/24M1658218}.

\bibitem{McDiarmid}
C.~McDiarmid.
\newblock On the method of bounded differences.
\newblock In J.~Siemons, editor, {\em Surveys in Combinatorics, 1989}, volume
  141 of {\em London Math. Soc. Lecture Note Ser.}, pages 148--188. Cambridge
  Univ. Press, Cambridge, 1989.
\newblock
  \href{https://doi.org/10.1017/CBO9781107359949.008}{doi:10.1017/CBO9781107359949.008}.

\bibitem{MubayiTuranExtension}
D.~Mubayi.
\newblock A hypergraph extension of {Tur{\'a}n}'s theorem.
\newblock {\em J. Combin. Theory Ser. B}, 96(1):122--134, 2006.
\newblock
  \href{https://doi.org/10.1016/j.jctb.2005.06.013}{doi:10.1016/j.jctb.2005.06.013}.

\bibitem{Ramsey}
F.~P. Ramsey.
\newblock On a problem of formal logic.
\newblock {\em Proc. London Math. Soc. (2)}, 30(1):264--286, 1930.
\newblock
  \href{https://doi.org/10.1112/plms/s2-30.1.264}{doi:10.1112/plms/s2-30.1.264}.

\bibitem{ReiherSurvey}
C.~Reiher.
\newblock Extremal problems in uniformly dense hypergraphs.
\newblock {\em European J. Combin.}, 88:103117, 2020.
\newblock
  \href{https://doi.org/10.1016/j.ejc.2020.103117}{doi:10.1016/j.ejc.2020.103117}.

\bibitem{RRStetrahedra}
C.~Reiher, V.~R{\"o}dl, and M.~Schacht.
\newblock Embedding tetrahedra into quasirandom hypergraphs.
\newblock {\em J. Combin. Theory Ser. B}, 121:229--247, 2016.
\newblock
  \href{https://doi.org/10.1016/j.jctb.2016.06.008}{doi:10.1016/j.jctb.2016.06.008}.

\bibitem{RRSzero}
C.~Reiher, V.~R{\"o}dl, and M.~Schacht.
\newblock Hypergraphs with vanishing {Tur{\'a}n} density in uniformly dense
  hypergraphs.
\newblock {\em J. Lond. Math. Soc. (2)}, 97(1):77--97, 2018.
\newblock \href{https://doi.org/10.1112/jlms.12095}{doi:10.1112/jlms.12095}.

\bibitem{RRS}
C.~Reiher, V.~R{\"o}dl, and M.~Schacht.
\newblock On a generalisation of {Mantel}'s theorem to uniformly dense
  hypergraphs.
\newblock {\em Int. Math. Res. Not. IMRN}, 2018(16):4899--4941, 2018.
\newblock \href{https://doi.org/10.1093/imrn/rnx017}{doi:10.1093/imrn/rnx017}.

\bibitem{RRSweak}
C.~Reiher, V.~R{\"o}dl, and M.~Schacht.
\newblock On a {Tur{\'a}n} problem in weakly quasirandom {$3$}-uniform
  hypergraphs.
\newblock {\em J. Eur. Math. Soc. (JEMS)}, 20(5):1139--1159, 2018.
\newblock \href{https://doi.org/10.4171/JEMS/784}{doi:10.4171/JEMS/784}.

\bibitem{RRStwopath}
C.~Reiher, V.~R{\"o}dl, and M.~Schacht.
\newblock Some remarks on {$\pi_{\wedge}$}.
\newblock In S.~Butler, J.~Cooper, and G.~Hurlbert, editors, {\em Connections
  in Discrete Mathematics: A Celebration of the Work of Ron Graham}, pages
  214--239. Cambridge Univ. Press, Cambridge, 2018.
\newblock
  \href{https://doi.org/10.1017/9781316650295.014}{doi:10.1017/9781316650295.014}.

\bibitem{RScounting}
V.~R{\"o}dl and M.~Schacht.
\newblock Regular partitions of hypergraphs: counting lemmas.
\newblock {\em Combin. Probab. Comput.}, 16(6):887--901, 2007.
\newblock
  \href{https://doi.org/10.1017/S0963548307008565}{doi:10.1017/S0963548307008565}.

\bibitem{RSregular}
V.~R{\"o}dl and M.~Schacht.
\newblock Regular partitions of hypergraphs: regularity lemmas.
\newblock {\em Combin. Probab. Comput.}, 16(6):833--885, 2007.
\newblock
  \href{https://doi.org/10.1017/S0963548307008553}{doi:10.1017/S0963548307008553}.

\bibitem{RSremoval}
V.~R{\"o}dl and M.~Schacht.
\newblock Generalizations of the removal lemma.
\newblock {\em Combinatorica}, 29(4):467--501, 2009.
\newblock
  \href{https://doi.org/10.1007/s00493-009-2320-x}{doi:10.1007/s00493-009-2320-x}.

\bibitem{Turan1941}
P.~Tur{\'a}n.
\newblock Egy gr{\'a}felm{\'e}leti sz{\'e}ls{\H{o}}{\'e}rt{\'e}kfeladatr{\'o}l.
\newblock {\em Mat. Fiz. Lapok}, 48:436--452, 1941.
\newblock In Hungarian.

\bibitem{WigdersonSidorenko}
Y.~Wigderson.
\newblock Complete {$r$}-partite {$r$}-graphs are {Sidorenko}: A brief
  exposition.
\newblock Expository note.
\newblock
  \href{https://ywigderson.math.ethz.ch/math/static/HypergraphSidorenko.pdf}{Available
  online}.

\end{thebibliography}
\endgroup

\end{document}